\documentclass[12pt,a4paper,fleqn]{article}
\usepackage{a4wide,amsfonts,amsmath,latexsym,amssymb,euscript,graphicx,units,mathrsfs}

\usepackage[english]{babel}
\usepackage{graphicx}
\usepackage{color}
\usepackage{amssymb}
\usepackage{amssymb}
\usepackage[T1]{fontenc}
\usepackage{latexsym}
\usepackage{xypic}
\usepackage{eufrak}
\usepackage{euscript}
\usepackage{amsfonts,amsmath}
\usepackage{verbatim}
\usepackage{fancyhdr}
\usepackage{mathrsfs}
\usepackage{units}

\newtheorem{prop}{Proposition}[section]

\newtheorem{lemme}[prop]{Lemma}
\newtheorem{lemma}[prop]{Lemma}

\newtheorem{remark}[prop]{Remark}

\newtheorem{theorem}[prop]{Theorem}

\renewcommand{\geq}{\geqslant}
\def\leq{\leqslant}

\newcommand{\N}{\mathbb{N}}
\newcommand{\Z}{\mathbb{Z}}

\newcommand{\R}{\mathbb{R}}

\newcommand{\floor}[1]{\left\lfloor#1\right\rfloor}

\newcommand{\bea}{\begin{equation}}
\newcommand{\eea}{\end{equation}}
\newcommand{\beas}{\begin{eqnarray*}}
\newcommand{\eeas}{\end{eqnarray*}}

\def\HH{\EuFrak H}

\def\e{\varepsilon}

\def\1{{\mathbf{1}}}

\def\1{{\mathbf{1}}}
\def\0.5{{\frac{1}{2}}}

\newcommand{\qed}{\nopagebreak\hspace*{\fill}
{\vrule width6pt height6ptdepth0pt}\par}

\begin{document}

\begin{center}
{\bf\Large Symmetric  weighted odd-power variations of  fractional Brownian motion and applications }
\end{center}
\bigskip
\begin{center}
{\bf David Nualart\footnote{The University of Kansas,  Lawrence, Kansas 66045, USA. Email: nualart@ku.edu.
} \:\:\: Raghid Zeineddine\footnote{Freiburg Institute For Advanced Studies and University of Strasbourg Institute for Advanced Studies (FRIAS-USIAS), Freiburg, Germany. Email: raghid.z@hotmail.com.} }
\end{center}

\bigskip
\begin{abstract}
We  prove a non-central limit theorem for the  symmetric weighted odd-power variations of the fractional Brownian motion  with Hurst parameter $H< 1/2$. As applications, we study the asymptotic behavior of  the trapezoidal weighted odd-power variations of the fractional Brownian motion and the fractional Brownian motion in Brownian time  $Z_t:=  X_{Y_t}, \:\: t \geq 0$, where $X$ is a fractional Brownian motion  and $Y$ is an independent Brownian motion. 
\end{abstract}

\section{Introduction}

Let $X=(X_t)_{ t \geq 0}$ be a fractional Brownian motion (fBm) with Hurst parameter $H\in (0,1/2)$. 
The purpose of this paper is to prove a non-central limit theorem for  symmetric weighted odd-power variations of $X$ and derive some applications.

For any integers  $n\geq 1$ and   $j\geq 0$ we will make use of the notation
 $\Delta_{j, n} X:=  X_{(j+1)2^{-n}}-X_{j2^{-n}}$ and   $\beta_{j,n}:=\frac 12(X_{j2^{-n}}+X_{(j+1) 2^{-n}})$. 
The main result of the paper is the following theorem.

\begin{theorem}\label{intro-main}  Let $X$ be a fBm with Hurst parameter $H< 1/2$. Fix an integer $r\geq 1$. Assume that $f\in C^{2r-1} (\mathbb{R})$.   Then, as  $n \to \infty$, we have
\begin{equation}
\left(2^{-n/2}\sum_{j=0}^{\lfloor 2^n t\rfloor -1 }f(\beta_{j,n})\big(2^{nH}(\Delta_{j,n} X)\big)^{2r-1}  \right)_{t\geq 0}\overset{\rm Law}{\longrightarrow}  \left( \sigma_r \int_0^t f(X_s)dW_s \right)_{t\geq 0} , 
\end{equation}
where   $W$ is a standard Brownian motion independent of $X$, $\sigma_r$ is the constant  given by
\begin{equation}   \label{sigma}
\sigma_r^2= E [ X_1^{4r-2} ] + 2 \sum_{j=1 } ^\infty E [ (X_1(X_{1+j}-X_j) ) ^{2r-1} ],
\end{equation}
and the convergence holds in the  Skorohod space $D([0, \infty))$.
\end{theorem}

 The proof of this result is based on the methodology of big blocks-small blocks, used, for instance, in  \cite{CNP,CNW} and the following stable convergence of odd-power variations  of the fBm  
   \begin{equation}  \label{er1}
 \left(    2^{-n/2} \sum_{j=0}^{\lfloor 2^nt \rfloor -1} \big( 2^{nH}\Delta _{j,n} X\big)^{2r-1} , X_t \right)_{t\geq 0}\,\mathrel{\mathop{\longrightarrow}^{\mathrm{Law}}_{n\rightarrow\infty}} \, \left( \sigma_r  W_t, X_t \right)_{t\geq 0},
 \end{equation}
 where $\sigma _r$ is defined in (\ref{sigma}) and  in the right-hand side of (\ref{er1}), the process $W$ is a Brownian motion  independent of $X$. The proof of the convergence   (\ref{er1})  for a fixed $t$ follows from the  Breuer-Major Theorem
 (we refer to  \cite[Chapter 7]{NP2}   and \cite{CNW} for a proof of this result  based on the Fourth Moment theorem). 
 
 \medskip
   A rather complete analysis of the asymptotic behavior of weighted power variations of the fBm was developed in \cite[Corollary 3]{NoNuTu}. However, the case  of symmetric weighted power  variations  was not considered in this paper.  On the other hand,  motivated by applications to the asymptotic behavior of symmetric Riemann sums for critical values of the Hurst parameter, 
Theorem  1.1 was proved in \cite[Proposition 3.1]{BNN}  when $H=\frac 1{ 4r-2}$ for a function of the form $f^{(2r-1)}$ and assuming that
$f\in C^{20r-15}(\R)$  is such that $f$ and its derivatives up to the order $20r-15$ have moderate growth.  The proof given here,  inspired by the recent work of Harnett, Jaramillo and Nualart \cite{HJN},  allows less derivatives and no growth condition.

\medskip
In the second part of the paper we present two applications of Theorem 1.1
First, we deduce the  following  convergence in law of the trapezoidal weighted odd-power variations of the fBm  with Hurst parameter $H < 1/2$.  
\begin{prop} \label{Prop}
  Let $X$ be a fBm with Hurst parameter $H < 1/2$. Fix  an integer $r\geq 1$. Then, if $f\in C^M(\mathbb{R})$, where $M >2r-2 + \frac 1{2H} $, 
  as $n\rightarrow \infty$, we have
\[
 \left(2^{-n/2}\sum_{j=0}^{\lfloor 2^{n} t\rfloor -1}\frac12(f(X_{j2^{-n}})+ f(X_{(j+1)2^{-n}}))\big(2^{nH} \Delta_{j,n}X\big)^{2r-1} \right)_{t\geq 0}
 {\overset{Law}{\longrightarrow}}  \left(\sigma_r \int_0^t f(X_s)dW_s \right)_{t\geq 0},
\]
in the Skorohod space $D([0,\infty))$, where $W$ is a Brownian motion independent of $X$.

\end{prop}
 In the particular case $r=2$ and $H=1/6$, this result  has been proved in \cite{NRS} with  longer arguments and using in a methodology introduced in \cite{NoNu}. The limit in this case, that is $\sigma_2 \int_0^t f(X_s)dW_s$, is the correction term in the  It\^o-type formula in law proved in \cite{NRS}.

The asymptotic behavior of  weighted odd-power variations of fBm with Hurst parameter $H<1/2$ has been already studied (see \cite{NoNuTu} and the references therein). More precisely, it is proved that for $H<1/2$, for any integer $r \geq 2$, and for a sufficiently smooth function $f$, we have
\[
2^{nH-n}\sum_{j=0}^{\lfloor 2^{n} t\rfloor -1}f(X_{j2^{-n}})\big(2^{nH} \Delta_{j,n}X\big)^{2r-1} 
 \underset{n\to \infty}{\overset{L^2}{\longrightarrow}} -\frac{\mu_{2r}}{2} \int_0^t f'(X_s)ds, 
\]  
where $\mu_{2r}:=E[N^{2r}]$ with $N \sim \mathscr{N}(0,1)$. By similar arguments, one can show that
\[
2^{nH-n}\sum_{j=0}^{\lfloor 2^{n} t\rfloor -1}f(X_{(j+1)2^{-n}})\big(2^{nH} \Delta_{j,n}X\big)^{2r-1} 
 \underset{n\to \infty}{\overset{L^2}{\longrightarrow}} \frac{\mu_{2r}}{2} \int_0^t f'(X_s)ds, 
\] 
which implies that 
 \begin{equation} \label{ecua1}
2^{nH-n}\sum_{j=0}^{\lfloor 2^{n} t\rfloor -1}\frac12(f(X_{j2^{-n}})+ f(X_{(j+1)2^{-n}}))\big(2^{nH} \Delta_{j,n}X\big)^{2r-1} 
 \underset{n\to \infty}{\overset{L^2}{\longrightarrow}} 0. 
\end{equation}
Thus, a natural question is to know  whether it is possible to replace the normalization $2^{nH-n}$ by another one in order to get a non-degenerate limit in the  convergence (\ref{ecua1})?   Proposition 1.2 gives us the answer to this question.
 
 \medskip
Our second application of Theorem 1.1 deals with  the asymptotic behavior of  weighted odd-power variations of the so-called {\it fractional Brownian motion in Brownian time} (fBmBt in short) when $H<1/2$. The fBmBt  is defined as
\[
Z_t=X_{Y_t}, \quad t\geq 0,
\]
where $X$ is a two-sided fractional Brownian motion, with Hurst parameter $H \in (0,1)$, and $Y$ is a standard (one-sided) Brownian
motion independent of $X$.  The process $Z$ is  self-similar of order $H/2$, it has  stationary increments but it is not Gaussian. In  the case $H=\frac12$,  where $X$ is a {\it standard} Brownian motion, one  recovers the celebrated \textit{iterated Brownian motion} (iBm). This terminology was coined by Burdzy in 1993 (see \cite{burdzy ibm}), but the idea of  considering the  iBm is actually older than that. Indeed,
Funaki \cite{funaki79} discovered  in 1979 that iBm may be used to represent the solution of 
the following parabolic partial differential equation:
\[
\frac{\partial u}{\partial t} = \frac18 \bigg( \frac{\partial u}{\partial x} \bigg)^4, \: (t,x) \in (0, \infty)\times \R.
\]
We refer the interested reader to the research works of Nane (see, e.g., \cite{erkane-nane} and the references therein)
for many other interesting relationships between iterated processes and partial differential equations.

In 1998, Burdzy and Khoshnevisan \cite{3} showed that iBm can be somehow considered as the canonical motion in an independent Brownian fissure. As such, iBm reveals to be a suitable candidate to model a diffusion in a Brownian crack. To support their claim, they have shown that the two components of a reflected two-dimensional   Brownian motion in a Wiener sausage of width $\epsilon>0$ converge to the usual Brownian motion and iterated Brownian motion, respectively, when $\epsilon$ tends to zero.

Let us go back to the second application of Theorem 1.1, we have the following theorem on the convergence in law of modified weighted
odd-power variations of the fBmBt. 

\begin{theorem}\label{application theorem}
Suppose that $H<  \frac12$ and fix an integer $r\geq 1$. Let $f\in C^M(\mathbb{R})$, where $M >2r-2 + \frac 1{2H} $.  Then, we have
\begin{eqnarray}
&&\bigg( 2^{-\frac{n}{4}}\sum_{k=0}^{\lfloor 2^n t \rfloor -1} \frac12\big(f(Z_{T_{k,n}})+f(Z_{T_{k+1,n}})\big)\big(2^{\frac{nH}{2}}(Z_{T_{k+1,n}}-Z_{T_{k,n}})\big)^{2r-1}  \bigg)_{t \geq 0} \notag\\
&&\underset{n\to \infty}{\overset{\rm Law}{\longrightarrow}} \bigg( \sigma_r \int_0^{Y_t} f(X_s)dW_s \bigg)_{t\geq 0}, \label{fdd}
\end{eqnarray}
in the Skorohod space $D([0, \infty))$, where  for $u \in \R$, $\int_0^{u} f(X_s)dW_s$ is the Wiener-It\^o integral of $f(X)$ with respect to $W$ defined in (\ref{integrale}) and  $\{T_{k,n}:\,1\leq k \leq 2^nt\}$ is a collection of stopping times defined in (\ref{T-N}) that approximates the
common dyadic partition $\{k2^{-n}:\,1\leq k \leq 2^nt\}$ of order $n$ of the time interval $[0,t]$.
\end{theorem} 
 
 Theorem \ref{application theorem} completes the study of the asymptotic behavior of the modified weighted odd-power variations of   the fBmBt in \cite{RZ4},   where the case $H \leq 1/6$ was missing. In addition, in Theorem \ref{application theorem} we have convergence in the Skorohod space  $D([0, \infty))$,  whereas in \cite{RZ4} we only proved the convergence of the finite dimensional distributions.

We remark that in many papers (see, for instance \cite{BNN})  the authors use the uniform partition, but in this paper we work with dyadic partitions.
Actually,  Theorem \ref{intro-main}  and Proposition 1.2 hold  also with the uniform partition. However,  the dyadic partition  plays a crucial  role in  Theorem \ref{application theorem}.

 The paper is organized as follows. In Section 2 we give some elements of Malliavin calculus and some preliminary results. In Section 3, we prove Theorem \ref{intro-main} and finally in Section 4 we prove Proposition \ref{Prop} and Theorem \ref{application theorem}.

\section{Elements of Malliavin calculus}

 In this section, we gather some elements of Malliavin calculus we shall need in the sequel. The reader in referred to \cite{Nualart, NP2} for details and any unexplained result.

Suppose that  $X = (X_{t})_{ t \in \R}$ a two-sided fractional Brownian motion with Hurst parameter $ H \in (0,1).$ That is, $X$ is a zero mean Gaussian process, defined on a complete probability space $(\Omega, \mathscr{A}, P)$, with  covariance function, \[ C_{H}(t,s) = E(X_{t}X_{s})=\frac{1}{2}(|s|^{2H} + |t|^{2H} -|t-s|^{2H}),\text{\:\:\:} s,t \in \R.\]
   We suppose that $\mathscr{A}$ is the $\sigma$-field generated by $X$. For all $n \in \N^*$, we let $\mathscr{E}_n$ be the set of step functions on $[-n,n]$, and $\displaystyle{\mathscr{E}:= \cup_n \mathscr{E}_n}$. Set $\varepsilon_t = \textbf{1}_{[0,t]}$ (resp. $\textbf{1}_{[t,0]}$) if $t \geq 0$ (resp. $t < 0$). Let $\mathfrak{H}$ be the Hilbert space defined as the closure of $\mathscr{E}$ with respect to the inner product
   \begin{equation}
    \langle \varepsilon_t, \varepsilon_s \rangle_{\mathfrak{H}} = C_{H}(t,s),\quad s,t \in \R. \label{inner product}
   \end{equation}
     The mapping $\varepsilon_t \mapsto X_{t}$ can be extended to an isometry between $\mathfrak{H}$ and the Gaussian space $\mathbb{H}_{1}$ associated with $X$. We will denote this isometry by $\varphi \mapsto X(\varphi).$

     Let $\mathscr{F}$ be the set of all smooth cylindrical random variables  of the form
   \[ 
   F = \phi (X_{t_{1}},\dots,X_{t_{l}}),
   \] where $l \in \N^*$, $\phi : \mathbb{R}^{l}\rightarrow \mathbb{R}$ is a $C^{\infty}$-function such that $f$ and all  its partial derivatives have at most polynomial growth, and  $ t_{1} <  \cdots <t_{l}$ are some real numbers. The derivative of $F$ with respect to $X$ is the element of $L^{2}(\Omega; \mathfrak{H})$ defined by 
   \[ 
   D_{s}F = \sum_{i=1}^{l}\frac{\partial\phi}{\partial x_{i}}(X_{t_{1}}, \dots ,X_{t_{l}})\varepsilon_{t_{i}}(s), \text{\: \: \:} s \in \R.
   \]
   In particular $D_{s}X_{t} = \varepsilon_t(s)$. For any integer $k \geq 1$, we denote by $\mathbb{D}^{k,2}$ the closure of $\mathscr{F}$ with respect to the norm
   \[ \|F\|_{k,2}^{2} = E(F^{2}) + \sum_{j=1}^{k} E[ \|D^{j}F\|_{\mathfrak{H}^{\otimes j}}^{2}].\]
   The Malliavin derivative $D$ satisfies the chain rule. If $\varphi : \mathbb{R}^{n} \rightarrow \mathbb{R}$ is $C_{b}^{1}$ and if $F_1,\ldots,F_n$ are in $\mathbb{D}^{1,2}$, then $\varphi(F_{1},\dots,F_{n}) \in \mathbb{D}^{1,2}$ and we have
   \[ D\varphi(F_{1},\dots,F_{n}) = \sum_{i=1}^{n} \frac{\partial \varphi}{\partial x_{i}}(F_{1},\dots,F_{n})DF_{i}.\]

   We denote by $\delta$  the adjoint of the derivative operator $D$,   also called the divergence operator. A random element $u \in L^{2}(\Omega; \mathfrak{H})$ belongs to the domain of the divergence operator $\delta$, denoted Dom$(\delta)$, if and only if it satisfies
   \[ |E\langle DF,u\rangle_{\mathfrak{H}}|\leq c_{u}\sqrt{E(F^{2})} \text{\: for  any\:} F\in \mathscr{F}.\]
   If $u \in$ Dom$(\delta)$, then $\delta(u)$ is defined by the duality relationship
\begin{equation}
E \big( F\delta(u)\big) = E \big( \langle DF,u\rangle_{\mathfrak{H}}\big), \label{adjoint1}
\end{equation}
for every $F \in \mathbb{D}^{1,2}$.

For every $n\geq 1$, let $\mathbb{H}_{n}$ be the $n$th Wiener chaos of $X$, that is, the closed linear subspace of $ L^{2}(\Omega, \mathscr{A},P)$ generated by the random variables $\lbrace H_{n}(X(h)), h \in \mathfrak{H}, \|h\|_{\mathfrak{H}}=1 \rbrace,$ where $H_{n}$ is the $n$th Hermite polynomial. Recall that  $H_0=0$, $H_p(x)= (-1)^p \exp(\frac{x^2}{2})\frac{d^p}{dx^p}\exp(-\frac{x^2}{2})$ for $p\geq 1$.
 The mapping
\begin{equation}
I_{n}(h^{\otimes n}) := H_{n}(X(h))\label{linear-isometry}
\end{equation}
 provides a linear isometry between the symmetric tensor product $\mathfrak{H}^{\odot n}$ and $\mathbb{H}_{n}$.  The  relation (\ref{adjoint1}) extends to the multiple Skorohod integral $\delta^q$ $(q \geq 1)$, and we have 
\begin{eqnarray}\label{adjoint}
  E \big( F\delta^q(u)\big) = E \big( \langle D^{q}F,u\rangle_{\mathfrak{H}^{\otimes q}}\big),
  \end{eqnarray}
  for any element $ u$ in the domain of $\delta^q$, denoted Dom$(\delta^q)$, and any random variable $F \in \mathbb{D}^{q,2}.$ Moreover, $\delta^q(u)=I_{q}(u)$ for any $u \in \mathfrak{H}^{\odot q}$.

For any Hilbert space $V$, we denote $\mathbb{D}^{k,p}(V)$ the corresponding Sobolev space of $V$-valued random variables (see \cite[page 31]{Nualart}). The operator $\delta^q$ is continuous from $\mathbb{D}^{k,p}(\mathfrak{H}^{\otimes q})$ to $\mathbb{D}^{k-q,p}$, for any $p>1$ and every integers $k \geq q \geq 1$, that is, we have
\[
\|\delta^q(u)\|_{\mathbb{D}^{k-q,p}} \leq C_{k,p} \|u\|_{\mathbb{D}^{k,p}(\mathfrak{H}^{\otimes q})},
\]
for all $u \in \mathbb{D}^{k,p}(\mathfrak{H}^{\otimes q})$ and some constant $C_{k,p}>0$. These estimates are consequences of Meyer inequalities (see \cite[Proposition 1.5.7]{Nualart}). We need the following result (see \cite[Lemma 2.1]{NoNu}) on the Malliavin calculus with respect to $X$. 
 \begin{lemma}\label{lem:Fdelta}
Let $q\geq1$ be an integer. Suppose that $F\in\mathbb{D}^{q,2}$, and let $u$ be a symmetric element in ${\rm Dom} \,\delta^{q}$. Assume that, for any $0\leq r+j\leq q$, $\left  \langle D^{r}F,\delta^{j}(u)  \right\rangle_{\HH^{\otimes r}}\in L^{2}(\Omega;\HH^{\otimes q-r-j})$. Then, for any $r=0,\dots, q-1$, $\left \langle D^{r}F,u \right\rangle_{\HH^{\otimes r}}$ belongs to the domain of $\delta^{q-r}$ and we have 
\begin{align*} 
F\delta^{q}(u)
  &=\sum_{r=0}^{q}\binom{q}{r}\delta^{q-r}(\left \langle D^{r}F,u \right\rangle_{\HH^{\otimes r}}).
\end{align*}
\end{lemma}

\medskip

 Let $\lbrace e_{k}, k \geq 1\rbrace $ be a complete orthonormal system in $\mathfrak{H}.$ Given $f \in \mathfrak{H}^{\odot n}$ and $g \in \mathfrak{H}^{\odot m},$ for every $r= 0,\dots,n\wedge m,$ the contraction of $f$ and $g$ of order $r$ is the element of $ \mathfrak{H}^{\otimes(n+m-2r)}$ defined by
  \[ 
  f\otimes_{r} g = \sum_{k_{1},\dots,k_{r} =1}^{\infty} \langle f, e_{k_{1}}\otimes \cdots \otimes e_{k_{r}}\rangle_{\mathfrak{H}^{\otimes r}}\otimes \langle g,e_{k_{1}}\otimes \cdots\otimes e_{k_{r}}\rangle_{\mathfrak{H}^{\otimes r}}.
  \]

\subsection{Preliminary results}
We will make use of the following notation:
\[
\partial_{j2^{-n}}  =\mathbf{1}_{[j2^{-n}, (j+1) 2^{-n}]}, \quad \e_t =\mathbf{1}_{[0,t]} ,\quad
\widetilde{\e} _{j2^{-n}} = \frac 12 \left( \e_{j2^{-n}} + \e_{(j+1) 2^{-n}} \right).
\]
We need the following preliminary results.
\begin{lemma}
We fix two integers $n>m\geq 2$,  and for any $j\geq 0$,  define
  $k:=k(j)= \sup \{i\geq 0:  i2^{-m}  \leq   j2^{-n}\}$.
The following inequality holds true for some constant $C_T$ depending only on $T$:
\bea \label{2.6}
\sum_{j=0}^{ \lfloor 2^n T\rfloor -1} \big| \big\langle\partial_{j2^{-n}}, \widetilde{\e}_{k(j) 2^{-m}}\big\rangle_{\mathfrak{H} }\big|    \leq  C_T
 2^{m(1-2H)}.\\  
\eea
 \end{lemma}
 
 \noindent
{\it Proof.}
See Lemma 2.2, inequality (2.11), in the paper by Binotto Nourdin and Nualart \cite{BNN}. In this paper the inequality is proved for
$\e_{k(j) 2^{-m}}$ but the case   $\widetilde{\e}_{k(j) 2^{-m}}$ can be proved by the same arguments.
 \qed
 
\begin{lemma} Let $0\leq s<t$. Then
 \bea
 \label{2.7}
     \sum_{j=\lfloor 2^n s\rfloor }^{\lfloor 2^n t\rfloor -1} \big| \big\langle\partial_{j2^{-n} },\widetilde{\varepsilon}_{j2^{-n}} \big\rangle_\mathfrak{H} \big| =  
     \frac 12 2^{-2nH} \left( \floor{2^n t}-\floor{2^n s} \right)^{2H}.
\eea
\end{lemma}

 \noindent
{\it Proof.}
We can write
\beas
       \sum_{j=\lfloor 2^n s\rfloor }^{\lfloor 2^n t\rfloor -1} \big| \big\langle\partial_{j2^{-n} },\widetilde{\varepsilon}_{j2^{-n}} \big\rangle_\mathfrak{H} \big|  &=&\frac12       \sum_{j=\lfloor 2^n s\rfloor }^{\lfloor 2^n t\rfloor -1} \Big| \mathbb{E}\big[ \big(X_{(j+1)2^{-n}}- X_{j2^{-n}}\big) \big(X_{(j+1)2^{-n}}+ X_{j2^{-n}}\big) \big] \Big|
\\
&=&   2^{-2nH}\frac 12       \sum_{j=\lfloor 2^n s\rfloor }^{\lfloor 2^n t\rfloor -1} \left[ (j+1)^{2H} -  j^{2H}\right],
\eeas
which gives the desired result.
\qed

\section{Proof of Theorem \ref{intro-main}}
In this section we provide the proof of Theorem  \ref{intro-main}.
We will make use of the following notation:
\bea  \label{jj2}
\Phi_n(t)= 2^{-n/2}\sum_{j=0}^{\lfloor 2^n t\rfloor -1}f(\beta_{j,n})\big(2^{nH}   \Delta_{j,n} X \big)^{2r-1}
\eea
and  $Z_t= \sigma_r \int_0^t f(X_s)dW_s$, where we recall that $W$ is a Brownian motion independent of $X$.
In order to prove Theorem 1.1, we need to show the following two results:

\medskip
\noindent
{\bf (A)}  {\it Convergence of the finite dimensional distributions:} Let   $0\leq t_{1}< \cdots <  t_{d}$ be fixed.   Then, we have
\[
 (\Phi_{n}(t_{1}),\dots, \Phi_{n}(t_{d}))     \stackrel{ Law} {\rightarrow}  (Z_{t_{1}},\dots, Z_{t_d}).
\]

\medskip
\noindent
{\bf (B)}  {\it Tightness:} The sequence $\Phi_{n}$ is tight in $D([0,\infty))$.  That is,
for every $\varepsilon, T>0$, there is a compact set $K\subset  D([0, T])$, such that 
\[
\sup_{n\geq 1}  P\left[\Phi_{n}\in K^{c}\right]<\varepsilon.
\]

 The proof of statements {\bf (A)} and {\bf (B)} will be done in several steps.

 \subsubsection*{Step 1: Reduction to compact support functions}

  As in \cite{HJN} in the proof of {\bf (A)} and {\bf (B)}  we can assume that $f$ has compact support. Indeed, fix $L\geq 1$ and let $f_L\in C^{2r-1}(\R)$ be a compactly supported function, such that $f_L(x) =f(x)$ for all $x\in [-L,L]$. Define
 \bea  \label{jj2}
\Phi_n ^L(t)= 2^{-n/2}\sum_{j=0}^{\lfloor 2^n t\rfloor -1}f_L(\beta_{j,n})\big(2^{nH}   \Delta_{j,n} X \big)^{2r-1}
\eea
and  $Z_t^L= \sigma_r \int_0^t f_L(X_s)dW_s$.  For  {\bf (B)}, we choose $L$ such that $P(\sup_{t\in [0,T]} |X_t| >L) <\frac {\varepsilon} 2$. Then, if
 $K_L \subset D([0,T])$ is a compact set that  $\sup_{n\geq 1}  P \left[\Phi^L_{n}\in K_L^{c}\right]<\frac \varepsilon 2$, we obtain
 \[
 P \left[\Phi_{n}\in K_L^{c}\right] \leq P\left[ \Phi^L_{n}\in K_L^{c}  , \sup_{t\in [0,T]} |X_t| \leq L \right] + P(\sup_{t\in [0,T]} |X_t| >L)  <\varepsilon.
 \]
With a similar argument, we can show that given a compactly supported function $\phi \in C(\R^d)$, the limit
\[
 \lim _{n\rightarrow \infty} E[ \phi(\Phi^L_{n}(t_{1}),\dots, \Phi^L_{n}(t_{d}))   - \phi(Z^L_{t_{1}},\dots, Z^L_{t_d})] =0
 \]
 implies the same limit with  $\Phi^L_{n}(t_{i})$ replaced by $ \Phi_{n}(t_{i})$ and $Z^L_{t_i}$ replaced by $Z_{t_i}$.

\subsubsection*{Step 2: Proof  {\bf (A)} assuming that $f$ has compact support}

  The proof is based on the small blocks-big blocks approach. Fix $m\leq n$  and for each $j \geq 0$ we write $k:= k(j) =
\sup\{i\geq 0: i2^{-m}  \leq j2^{-n} \}$, that is, $k(j) $ is the largest dyadic number in the $m$th generation which is less or equal than $j2^{-n}$.
Define   
\bea \label{jj3}
\widetilde{\Phi} _{n,m} (t) =2^{-n/2}  \sum_{j=0}^{\lfloor 2^n t\rfloor -1}  f(\beta_{k(j),m}) \big(2^{nH}   \Delta_{j,n} X \big)^{2r-1}.
\eea
This term can be decomposed as follows 
\begin{eqnarray*}
\widetilde{\Phi} _{n,m} (t) &=&2^{-n/2}  \sum_{k=0}^{\lfloor 2^m t\rfloor -1}  f(\beta_{k,m})\sum_{j=k2^{n-m}}^{(k+1 )2^{n-m}-1}\big(2^{nH}   \Delta_{j,n} X \big)^{2r-1} \\
&&
+2^{-n/2}  f(\beta_{\lfloor 2^mt \rfloor,m} ) \sum_{j=\lfloor 2^mt \rfloor 2^{n-m}}^{ \lfloor 2^n t \rfloor -1}\big (2^{nH}\Delta_{j,n}X\big) ^{2r-1}
\end{eqnarray*}
 The convergence (\ref{er1}) implies   that for any $ \mathscr{A}$-measurable and bounded random variable $\eta$, the random vector
 $(\widetilde{\Phi} _{n,m} (t_1), \dots, \widetilde{\Phi} _{n,m} (t_d) , \eta)$ converges in law, as $n$ tends to infinity, to the vector
 $(Y^1_m, \dots, Y^d_m, \eta)$, where
 \[
 Y^i_m =\sigma_r  \sum_{k=0}^{\lfloor 2^m t_i\rfloor -1} f(\beta_{k,m})\big(   \Delta_{k,m} W\big) + \sigma_r f(\beta_{\lfloor 2^m t_i \rfloor ,m})
 \big (W_{t_i}- W_{(\lfloor 2^m t_i \rfloor ) 2^{-m}})
\]
for $i=1,\dots, d$.   Clearly,  $Y^i_m$ converges in $L^2(\Omega)$, as $m$ tends to infinity, to $Z_{t_i}$ for $i=1,\dots, d$. 

Then, it suffices to show that
\begin{align}\label{lim:PhiPhitilde}
\lim_{m\rightarrow\infty}\limsup_{n\rightarrow\infty}\sum_{i=1}^{d}\| \Phi_{n}(t_{i})-\widetilde{\Phi}_{n,m}(t_{i})\|_{L^{2}(\Omega)}
  &=0.
\end{align}

Let  $c_{1,r},\dots,c_{r,r}$ will denote the coefficients of the Hermite expansion of $x^{2r-1}$, namely,  
\begin{align*}
x^{2r-1}
  &=\sum_{u=1}^{r}c_{u,r}H_{2(r-u)+1}(x).
\end{align*}
Then, we can write 
\begin{equation} \label{jj1}
  (2^{nH} \Delta_{j,n} X)^{2r-1} = \sum_{u=1}^{r} c_{u,r} H_{2(r-u)+1}\left(2^{nH} \Delta_{j,n} X \right) = \sum_{u=1}^r c_{u,r} 2^{nH(2(r-u)+1)}\delta^{2(r-u)+1}\left(\partial_{ j2^{-n}}^{\otimes 2(r-u)+1}  \right).
  \end{equation}
  Set $w: =w(u)=2(r-u)+1$. 
  Substituting (\ref{jj1}) into (\ref{jj2}), yields
\[
  \Phi_n(t_i) =   \sum_{u=1}^{r} c_{u,r}  \sum_{j=0}^{\lfloor 2^n t_i\rfloor-1 }f(\beta_{j,n}) 2^{-n/2 +w nH} \delta^w\left(\partial_{ j2^{-n}}^{\otimes w}  \right).
  \]
  On the other hand, (\ref{jj3}) can be also written as
  \[
  \widetilde{  \Phi}_{n,m} (t_i) =   \sum_{u=1}^{r} c_{u,r}  \sum_{j=0}^{  \lfloor 2^n t_i \rfloor-1 }   f(\beta_{k(j),m})
2^{-n/2 +w nH} \delta^w\left(\partial_{ j2^{-n}}^{\otimes w}  \right).
\]
With the help of Lemma \ref{lem:Fdelta} we can express these terms as  linear combinations of Skorohod integrals:
\[
  \Phi_n(t_i) =   \sum_{u=1}^{r} c_{u,r}    \sum_{\ell =0}^w {w\choose \ell}    \Theta^n_{u,\ell}(t_i)
  \]
  and
  \[
  \widetilde{  \Phi}_{n,m} (t_i) =   \sum_{u=1}^{r} c_{u,r}  \sum_{\ell =0}^w {w\choose \ell}    \widetilde{\Theta}^{n,m}_{u,\ell}(t_i)
\]
where
\[
  \Theta^n_{u,\ell}(t_i)=2^{-\frac n2 +w nH}      \sum_{j=0}^{\lfloor 2^n t_i\rfloor-1 }
  \delta^{w-\ell} \left(   f^\ell(\beta_{j,n}) \partial_{ j2^{-n}}^{\otimes (w-\ell)}   \langle \widetilde{\e}_{j2^{-n}} ,\partial_{j2^{-n}} \rangle 
  _\HH^\ell \right), 
  \]
  and
  \[
    \widetilde{\Theta}^{n,m}_{u,\ell}(t_i)  =  
  2^{-\frac n2 +w nH} 
   \sum_{j=0}^{ \lfloor 2^n t_i\rfloor -1 } 
     \delta^{w-\ell} \left(   f^\ell(\beta_{k(j),m}) \partial_{ j2^{-n}}^{\otimes (w-\ell)}   \langle \widetilde{\e}_{k(j)2^{-m}} ,\partial_{j2^{-n}} \rangle 
  _\HH^\ell \right).
\]
Then, it suffices to show that
\[
\lim_{m\rightarrow\infty}\limsup_{n\rightarrow\infty}\sum_{i=1}^{d}\|   \Theta^n_{u,\ell}(t_i) -     \widetilde{\Theta}^{n,m}_{u,\ell}(t_i) \|_{L^{2}(\Omega)} =0
\]
for all $1\leq u \leq r$ and $0\leq \ell \leq w$. We can decompose the difference  $\Theta^n_{u,\ell}(t_i) -     \widetilde{\Theta}^{n,m}_{u,\ell}(t_i)$ as follows

\[
\Theta^n_{u,\ell}(t_i) -     \widetilde{\Theta}^{n,m}_{u,\ell}(t_i) =
  2^{-\frac n2 +w nH} 
   \sum_{j=0}^{\lfloor 2^n t_i\rfloor-1 }      \delta^{w-\ell} \left(   F^{n,m}_{j,\ell} \partial_{ j2^{-n}}^{\otimes (w-\ell)}   \right)=:   T^{n,m}_{i,\ell} 
   \]
  where
\[
F^{n,m}_{j,\ell}=     f^\ell(\beta_{j,n})  \langle \widetilde{\e}_{j2^{-n}} ,\partial_{j2^{-n}} \rangle 
  _\HH^\ell  -   f^\ell(\beta_{k(j),m})  \langle \widetilde{\e}_{k(j)2^{-m}} ,\partial_{j2^{-n}} \rangle 
  _\HH^\ell. 
  \]
  
   By Meyer's inequality
  \beas
  \| T^{n,m}_{i,\ell}\|^2_2  &\leq& 
  C 2^{-n +2 nwH}   \sum_{h=0}^{w-\ell}  \left\| \sum_{j=0}^{\lfloor 2^n t_i\rfloor-1 } 
  D^hF^{n,m}_{k(j),j,\ell}  \otimes \partial_{ j2^{-n}}^{\otimes (w-\ell)}  \right\| ^2_{ L^2(\Omega; \HH^{\otimes(w-l+h)})} \\
&=& 
  C 2^{-n +2 nwH}   \sum_{h=0}^{w-\ell} 
\sum_{j_1,j_2=0}^{\lfloor 2^n t_i\rfloor-1 }     \mathbb{E} \left[
\langle D^h F^{n,m}_{j_1,\ell} , D^h F^{n,m}_{j_2,\ell}  \rangle_{\HH^{\otimes h}} \right]  \langle  \partial_{ j_12^{-n}}, \partial_{ j_22^{-n}} \rangle^{w-\ell}\\
&\leq&C 2^{-n +2 n\ell H}   \sum_{h=0}^{w-\ell} 
\sum_{j_1,j_2=0}^{\lfloor 2^n t_i\rfloor-1 }  
\|  D^hF^{n,m}_{j_1,\ell}\| _{L^2(\Omega;\HH^{\otimes h})}   \|  D^h F^{n,m}_{j_2,\ell}\| _{L^2(\Omega;\HH^{\otimes h})}    
| \rho_H(j_1-j_2)|^{w-\ell}.
\eeas
We will consider two different cases:

\noindent
{\it Case $w-\ell \geq 1$}:
  We can make the decomposition
\begin{eqnarray*}
F_{j,\ell}^{n,m}
  &=& f^{(\ell)}( \beta_{j,n} ) \langle \widetilde{\e}_{j2^{-n}}  ^{\otimes \ell}  -  \widetilde{\e}_{k(j)2^{-m}}  ^{\otimes \ell} , 
   \partial_{j2^{-n}}^{\otimes \ell} \rangle_{\HH^{\otimes \ell}}\\
	&& + \left(f^{(\ell)}(\beta_{j,n})-f^{(\ell)}  (\beta_{k(j),m} )\right)  \langle{\widetilde{\varepsilon}_{k(j)2^{-m}}},\partial_{j2^{-n}}\rangle_{\HH}^{\ell},
\end{eqnarray*}
and hence, we have
\begin{eqnarray*}
D^{h}F_{j,\ell}^{n,m}
	&=& f^{(\ell+h)}( \beta_{j,n})\widetilde{\varepsilon}_{j2^{-n}}^{\otimes h}\langle \widetilde{\varepsilon}_{j2^{-n}}^{\otimes \ell}-\widetilde{\e}  _{k(j)2^{-m}}^{\otimes \ell},\partial_{j2^{-n}}^{\otimes \ell} \rangle_{\HH^{\otimes \ell}}\\
	&&+f^{(\ell+h)}( \beta_{j,n})\left(\widetilde{\varepsilon}_{j2^{-n}}^{\otimes h}-\widetilde{\e}_{k(j)2^{-m}}^{\otimes h}\right) \langle \widetilde{\e}_{k(j)2^{-m}},\partial_{j2^{-n}} \rangle_{\HH}^{\ell}\\
	&&+\left(f^{(\ell+h)}( \beta_{j,n})-f^{(\ell+h)}( \beta_{k(j),m})\right)\widetilde{\e}_{k(j)2^{-m}}^{\otimes h}\langle  \widetilde{\e}_{k(j)2^{-m}},\partial_{j2^{-n}}\rangle_{\HH}^{\ell}.
\end{eqnarray*}
From the previous equality, and the compact support condition of $f$, we deduce that there exists a constant $C>0$, such that
\begin{eqnarray*}
  \left\| D^{h}F_{j,\ell}^{n,m} \right\|_{L^2( \Omega ; \HH^{\otimes h})} 
&\leq&  C  \left\| \widetilde{\varepsilon}_{j2^{-n}} \right\| _{\HH} ^h  \left\| \widetilde{\varepsilon}_{j2^{-n}}^{\otimes \ell}-\widetilde{\e}  _{k(j)2^{-m}}^{\otimes \ell} \right\|_{ \HH^{\otimes \ell}}  \left \| \partial_{j2^{-n}}^{\otimes \ell} \right\| _{\HH^{\otimes \ell}}\\
	&&+C \left \| \widetilde{\varepsilon}_{j2^{-n}}^{\otimes h}-\widetilde{\e}_{k(j)2^{-m}}^{\otimes h}\right\|  _{\HH^{\otimes h}}
	\left\|  \widetilde{\e}_{k(j)2^{-m}} \right\| _{\HH} ^\ell  \left\| \partial_{j2^{-n}}  \right\| _{\HH} ^\ell \\
	&&+\left \|  f^{(\ell+h)}( \beta_{j,n})-f^{(\ell+h)}( \beta_{k(j),m})\right\|   _2  \left\| \widetilde{\e}_{k(j)2^{-m}} \right\| _\HH^{ h+\ell}
	\left\|  \partial_{j2^{-n}}\right\|_{\HH}^{\ell}.
\end{eqnarray*}
Using Cauchy-Schwarz inequality we get, for  any natural number $p\geq 1$
\beas
\left\| \widetilde{\varepsilon}_{j2^{-n}}^{\otimes p}-\widetilde{\e}  _{k(j)2^{-m}}^{\otimes p} \right\|_{ \HH^{\otimes p}} 
&\leq&    \left\| \widetilde{\varepsilon}_{j2^{-n}}-\widetilde{\e}  _{k(j)2^{-m}}\right\|_{ \HH} 
\sum_{i=0}^{p-1}  \left\| \widetilde{\varepsilon}_{j2^{-n}}  \right\|_\HH ^i   \left\|    \widetilde{\e}  _{k(j)2^{-m}} \right\|_{\HH}  ^{  p-1-i}  \\
&\leq&  C  \left\| \widetilde{\varepsilon}_{j2^{-n}}-\widetilde{\e}  _{k(j)2^{-m}}\right\|_{ \HH} .
\eeas
Therefore
\beas
  \left\| D^{h}F_{j,\ell}^{n,m} \right\|_{L^2( \Omega ; \HH^{\otimes h})} 
 &\leq&    C \left\|  \partial_{j2^{-n}}\right\|_{\HH}^{\ell}
\left(  \left\| \widetilde{\varepsilon}_{j2^{-n}}-\widetilde{\e}  _{k(j)2^{-m}}\right\|_{ \HH}
+\left \|  f^{(\ell+h)}( \beta_{j,n})-f^{(\ell+h)}( \beta_{k(j),m})\right\|   _2   \right)\\
&\leq&
C 2^{-\ell nH} \left(   \sup_{|t-s| \leq 2^{-m}}  \left\|   X_t-X_s \right\|_2 + \left \|  f^{(\ell+h)}( \beta_{j,n})-f^{(\ell+h)}( \beta_{k(j),m})\right\|   _2  \right).
\eeas
Because $f^{(\ell+h)}$ is uniformly continuous, for any given $\e>0$ there exists $\delta>0$ such that $|x-y| <\delta$ implies $|f^{(\ell+h )} (x) - f^{(\ell +h)} (y) |<\e$. Therefore, we can write
\[
 \left \|  f^{(\ell+h)}( \beta_{j,n})-f^{(\ell+h)}( \beta_{k(j),m})\right\|   _2  
 \leq \e+  \frac  2 \delta  \| f^{(\ell +h)} \|_\infty   \| \beta_{j,n} -\beta_{k(j),m} \|_2 
 \]
 and this leads to the estimate
 \[
  \left\| D^{h}F_{j,\ell}^{n,m} \right\|_{L^2( \Omega ; \HH^{\otimes h})} 
  \leq C 2^{-\ell nH} \left(   \sup_{|t-s| \leq 2^{-m}}  \left\|   X_t-X_s \right\|_2 +\e \right),
  \]
  which implies
  \beas
  \| T^{n,m}_{i,\ell}\|^2_2
  &\leq&
    C   
  \left(   \sup_{|t-s| \leq 2^{-m}}      \sum_{i=0}^{w-\ell} \left\|   X_t-X_s \right\|_2 + \e \right)^2
 \sum_{j=0}^{\lfloor 2^n t_i\rfloor -1 }  
 |\rho_H(j)|^{w-\ell}.
\eeas
Then,  the series  $ \sum_{j=0}^{\infty}| \rho_H(j)|^{w-\ell}$ is convergent because $w-\ell \geq 1$ and $H<1/2$, and we obtain
 \[
 \lim_{m\rightarrow \infty}  \sup_n   \| T^{n,m}_{i,\ell}\|^2_2  =0,
 \]
 because $\e$ is arbitrary. 
 
 \noindent
 {\it Case  $\ell =w$}.  in this case  we have
 \beas
\| T^{n,m} _{i,w} \|^2_2&\leq &
   2^{-n +2w nH}  \left(
  \sum_{j=0}^{\lfloor   2^n  t_i\rfloor-1 } 
     \|  F^{n,m}_{j,w}   \|_2 \right)^2\\
   &\leq&  C     2^{-n +2w nH}   \left(
  \sum_{j=0}^{\lfloor    2^n   t_i\rfloor-1 } 
       |  \langle \widetilde{\e}_{j2^{-n}} ,\partial_{j2^{-n}} \rangle 
  _\HH^w  |   + | \langle \widetilde{\e}_{k(j)2^{-m}} ,\partial_{j2^{-n}} \rangle 
  _\HH^w    |
 \right)^2\\
  &\leq&  C     2^{n(2H-1) }   \left(
  \sum_{j=0}^{\lfloor    2^n  t_i\rfloor-1 } 
       |  \langle \widetilde{\e}_{j2^{-n}} ,\partial_{j2^{-n}} \rangle 
  _\HH  |   + | \langle \widetilde{\e}_{k(j)2^{-m}} ,\partial_{j2^{-n}} \rangle 
  _\HH    |
 \right)^2.
 \eeas
 Finally, using (\ref{2.6}) and (\ref{2.7}), we obtain
 \[
  \| T^{n,m} _{i,w} \|^2_2 \leq C  2^{n(2H-1) }  2^{2m(1-2H)},
  \]
  which implies
 \[
 \lim_{m\rightarrow \infty}  \limsup_{n\rightarrow \infty}   \| T^{n,m}_{i,w}\|^2_2  =0.
 \]

 \subsubsection*{Step 3:  Proof  {\bf (B)}  assuming that $f$ has compact support.}

We claim that for   every $0\leq s\leq t\leq T$, and $p>2$, there exists a constant $C>0$, such that
\bea \label{claim}
  E \left[ |\Phi_{n}(t)-\Phi_{n}(s)|^{p}\right]\leq C \left(\frac{\floor{2^n t}-\floor{2^n s}}{2^n}\right)^{\frac{p}{2}}
+ C\left(\frac{\floor{2^n t}-\floor{2^n s}}{2^n}\right)^{pH}.
\eea
Then, by the `Billingsley criterion' (see \cite[Theorem 13.5]{Billingsley}),  (\ref{claim}) implies the desired tightness property. 
From the computations in the proof of  {\bf (A)}, we need  to show that for any $1\leq u \leq r$ and for any $0\leq \ell \leq w$, where $w=2(r-u)+1$,
\bea \label{claim2}
\|  \Theta^n_{u,\ell} (t) - \Theta^n_{u,\ell} (s) \|_p \leq  C \left(\frac{\floor{2^n t}-\floor{2^n s}}{2^n}\right)^{\frac{1}{2}}
+ C\left(\frac{\floor{2^n t}-\floor{2^n s}}{2^n}\right)^{H}.
\eea
By Meyer's inequalities,
\beas
&& \|  \Theta^n_{u,\ell} (t) - \Theta^n_{u,\ell} (s) \|_p  \\
 &&\quad = 2^{-\frac n2 +w nH}  \left\|         \sum_{j=\lfloor 2^n s\rfloor}^{\lfloor 2^n t\rfloor-1 }
  \delta^{w-\ell} \left(   f^{(\ell)}(\beta_{j,n}) \partial_{ j2^{-n}}^{\otimes (w-\ell)}   \langle \widetilde{\e}_{j2^{-n}} ,\partial_{j2^{-n}} \rangle 
  _\HH^\ell \right) \right\|_p\\
 & & \quad \leq  C 2^{-\frac n2 +w nH}  \sum_{h=0} ^{w-\ell} \left\|         \sum_{j=\lfloor 2^n s\rfloor}^{\lfloor 2^n t\rfloor-1 }
     f^{(\ell+ h)} (\beta_{j,n})  \widetilde{\e} ^{\otimes h} _{j2^{-n}}  \otimes \partial_{ j2^{-n}}^{\otimes (w-\ell)}   \langle \widetilde{\e}_{j2^{-n}} ,\partial_{j2^{-n}} \rangle 
  _\HH^\ell  \right\|_{ L^p(\Omega; \HH^{\otimes (w-\ell +h)})} \\
  &&  \quad = C 2^{-\frac n2 +w nH}  \sum_{h=0} ^{w-\ell} \left\|     \left\|      \sum_{j=\lfloor 2^n s\rfloor}^{\lfloor 2^n t\rfloor-1 }
     f^{(\ell+ h)} (\beta_{j,n})  \widetilde{\e} ^{\otimes h} _{j2^{-n}}  \otimes \partial_{ j2^{-n}}^{\otimes (w-\ell)}   \langle \widetilde{\e}_{j2^{-n}} ,\partial_{j2^{-n}} \rangle 
  _\HH^\ell  \right\|^2_{\HH^{\otimes(w-\ell+h)}} \right\|_{\frac p2}^{\frac 12}.
 \eeas
 As a consequence,   since $f$ has compact support, applying Minkowski inequality, there is a constant $C$ such that
 \beas
 \|  \Theta^n_{u,\ell} (t) - \Theta^n_{u,\ell} (s) \|^2_p   &\leq &C 2^{-n+2w nH}  \sum_{h=0} ^{w-\ell} \bigg \|          \sum_{j,k=\lfloor 2^n s\rfloor}^{\lfloor 2^n t\rfloor-1 }
     f^{(\ell+ h)} (\beta_{j,n})   f^{(\ell+ i)} (\beta_{k,n})   \\
     &&  \times \left \langle  \widetilde{\e} _{j2^{-n}} ,  \widetilde{\e}_{k2^{-n}} \right \rangle^h
      \left \langle  \delta _{j2^{-n}} , \delta_{k2^{-n}} \right \rangle^{w-\ell}
    \left \langle \widetilde{\e}_{j2^{-n}} ,\partial_{j2^{-n}} \right  \rangle 
  _\HH^\ell   \left \langle \widetilde{\e}_{k2^{-n}} ,\partial_{k2^{-n}} \right  \rangle 
  _\HH^\ell
    \bigg\|_{\frac p2} \\
    & \leq &C 2^{-n +2\ell nH}          \sum_{j,k=\lfloor 2^n s\rfloor}^{\lfloor 2^n t\rfloor-1 }   
      \left| \rho_H(j-k) \right| ^{w-\ell}
 \left|   \left \langle \widetilde{\e}_{j2^{-n}} ,\partial_{j2^{-n}} \right  \rangle 
  _\HH \right| ^\ell   \left| \left \langle \widetilde{\e}_{k2^{-n}} ,\partial_{k2^{-n}} \right  \rangle 
  _\HH \right|^\ell.
 \eeas
 We will consider two different cases:
 
 \noindent
 {\it Case $ w-\ell \geq 1$}:
 In this case, we obtain
\beas
 \|  \Theta^n_{u,\ell} (t) - \Theta^n_{u,\ell} (s) \|_p^2    &\leq & C 2^{-n }        \sum_{j,k=\lfloor 2^n s\rfloor}^{\lfloor 2^n t\rfloor-1 }| \rho_H(j-k)|^{w-\ell} \\
 & \leq & C \frac{\floor{2^n t}-\floor{2^n s}}{2^n}  \sum_{h=\lfloor 2^n s\rfloor}^{\lfloor 2^n t\rfloor-1 } 
| \rho_H(h)|^{w-\ell} \\
 & \leq & C \frac{\floor{2^n t}-\floor{2^n s}}{2^n}.
 \eeas
because the series  $ \sum_{j=0}^{\infty}| \rho_H(j)|^{w-\ell}$ is convergent because $w-\ell \geq 1$ and $H<1/2$. This implies the inequality (\ref{claim2}) in this case.

 \noindent
 {\it Case $\ell =w$}:
 We have
 \beas
  \|  \Theta^n_{u,w} (t) - \Theta^n_{u,w} (s) \|_p^2   &\leq &C2^{-n+2w nH}  \left(
   \sum_{j=\lfloor 2^n s\rfloor}^{\lfloor 2^n t\rfloor-1 }   
 \left|   \left \langle \widetilde{\e}_{j2^{-n}} ,\partial_{j2^{-n}} \right  \rangle 
  _\HH \right| ^w   \right)^2 \\
  &\leq &C2^{-n+2 nH}  \left(
   \sum_{j=\lfloor 2^n s\rfloor}^{\lfloor 2^n t\rfloor-1 }   
 \left|   \left \langle \widetilde{\e}_{j2^{-n}} ,\partial_{j2^{-n}} \right  \rangle 
  _\HH \right|    \right)^2.
  \eeas
  Finally, applying (\ref{2.7}) and the fact that  $2^{-n+2 nH} \leq 1$, we obtain
  \[
   \|  \Theta^n_{u,w} (t) - \Theta^n_{u,w} (s) \|_p ^2  \leq C
    \left(  \frac {\floor{2^n t}-\floor{2^n s} } {2^n} \right)^{2H}.
   \]
   This completes the proof of part {\bf (B)}.

\section{Applications}

\subsection{The trapezoidal weighted odd-power variations of  fractional Brownian motion}
The trapezoidal weighted odd-power variations of  the fBm  is given in Proposition \ref{Prop}. We give its proof  below.

\medskip
\noindent
{\it Proof of Proposition 1.2.} 
  By a localization argument similar to that used in the  proof of Theorem  \ref{intro-main}, we can assume that $f$ has compact support.
 Choose an integer $N$ such that   $\frac 1{2H} -1 <N \leq M-(2r-1)$, which is possible because $M> 2r-2+\frac 1{2H}$.
Since $f\in C^M(\mathbb{R})$, by Taylor expansion, we have for all $x,y \in \R$ and $N \leq M-(2r-1)$, 
\begin{eqnarray*}
f(y)&=& f(\frac12(x+y)) + \frac12 f'(\frac12(x+y))(y-x) + \sum_{k=2}^{N}\frac{1}{2^k} \frac{1}{k!}f^{(k)}(\frac12(x+y))(y-x)^k  +R^{(1)}_N,\\
f(x)&=& f(\frac12(x+y)) + \frac12 f'(\frac12(x+y))(x-y) + \sum_{k=2}^{N}\frac{1}{2^k} \frac{1}{k!}f^{(k)}(\frac12(x+y))(x-y)^k +R^{(2)}_N,
\end{eqnarray*} 
where the residual terms $R^{(1)}_N$ and $R^{(2)}_N$ are bounded by  $C|y-x|^{N+1}$. We deduce that, for all integer $N \geq 1$,
\begin{equation}  \label{jj34}
\frac12(f(x)+f(y))= f(\frac12(x+y)) + \sum_{k=1}^{\lfloor\frac{N}{2}\rfloor}\frac{1}{2^k} \frac{1}{k!}f^{(2k)}(\frac12(x+y))(y-x)^{2k} +  R_N (x,y),
\end{equation}
where $R_N(x,y)  \leq  C|y-x|^{N+1} $.
Recall that  $\beta_{j,n}:=1/2(X_{j2^{-n}}+X_{(j+1)2^{-n}})$ and we also write $\Delta_{j,n} f(X):=  \frac12(f(X_{j2^{-n}})+ f(X_{(j+1)2^{-n}}))$.
Set 
\[
\Psi_n(t)=2^{-n/2}\sum_{j=0}^{\lfloor 2^{n} t\rfloor -1}\Delta_{j,n} f(X)\big(2^{nH} \Delta_{j,n} X\big)^{2r-1}
\]
and let $\Phi_n(t)$ be defined in (\ref{jj2}). Then, in view of Theorem 1.1, it suffices to show that the difference $\Psi_n- \Phi_n$  converges to zero in probability in  the Skorohod space as $n\rightarrow \infty$. Using the expansion (\ref{jj34}), we obtain
 
\begin{eqnarray*}
 \Psi_n(t)- \Phi_n(t) &=&
  2^{-n/2}\sum_{j=0}^{\lfloor 2^{n} t\rfloor -1}\big(\Delta_{j,n} f(X) - f(\beta_{j,n})\big)\big(2^{nH} \Delta_{j,n}X\big)^{2r-1}\\
&&= 2^{-n/2}2^{-2nHk} \sum_{k=1}^{\lfloor\frac{N}{2}\rfloor}\frac{1}{2^k} \frac{1}{k!}\sum_{j=0}^{\lfloor 2^{n} t\rfloor -1}f^{(2k)}(\beta_{j,n}) \big(2^{nH} \Delta_{j,n}X\big)^{2k+2r-1} \\
&& + 2^{-n/2}\sum_{j=0}^{\lfloor 2^{n} t\rfloor -1}R_N( X_{j2^{-n}},  X_{(j+1)2^{-n}})\big(2^{nH} \Delta_{j,n}X\big)^{2r-1} \\
&&=: A_n (t)+B_n(t)
\end{eqnarray*}
Thanks to Theorem \ref{intro-main},  and taking into account that  $f^{(2k)} \in C^{2k+2r-1}(\mathbb{R})$ for all $k \leq \lfloor N/2 \rfloor$
because $N +2r-1\leq M$, we deduce that $A_n(\cdot)$ converges to 0 in probability as $n \to \infty$ in $D([0,\infty))$. Therefore,  it is enough to prove the convergence in probability to 0 of $B_n(\cdot)$  in $D([0,\infty))$. This follows from the following estimates
\[
E\left[ \sup_{0 \leq t \leq T} |B_n(t)| \right]
\leq C 2^{-\frac{n}{2}}2^{-nH(N+1)}\sum_{j=0}^{\lfloor 2^{n} T\rfloor -1}  E [ |2^{nH} \Delta_{j,n} X |^{N+2r}] \leq C_T  2^{ \frac n2 -nH(N+1)},
\]
taking into account that $H(N+1)>\frac 12$. 
\qed

\subsection{The weighted power variations of fractional Brownian motion in Brownian time}
The so-called {\it fractional Brownian motion in Brownian time} (fBmBt in short) is defined as
\[
Z_t=X_{Y_t}, \quad t\geq 0,
\]
where $X$ is a two-sided fractional Brownian motion, with Hurst parameter $H \in (0,1)$, and $Y$ is a standard (one-sided) Brownian
motion independent of $X$.  The process $Z_t$  is not a Gaussian process and it is  self-similar (of order $H/2$) with stationary increments.
When $H=1/2$, one recovers the celebrated {\it iterated Brownian motion.}

 Let $f: \R\rightarrow \R$.   Then, for any $t\geq 0$ and any integer $p\geq 1$, the  weighted $p$-variation of $Z$ is defined as 
\begin{eqnarray*}
M_n^{(p)}(t)=\sum_{k=0}^{\lfloor 2^n t\rfloor-1}\frac12\big(f(Z_{k2^{-n}})+f(Z_{(k+1)2^{-n}})\big)( \Delta_{k,n}Z)^p. 
\end{eqnarray*} 
where, as before, $\Delta_{k,n} Z =Z_{(k+1)2^{-n}}- Z_{k2^{-n}}$.
 After proper normalization we may expect the  convergence (in some sense) to a non-degenerate limit (to be determined) of
\begin{equation}
N_n^{(p)}(t)= 2^{-n\kappa}\sum_{k=0}^{\lfloor 2^n t\rfloor-1}\frac12\big(f(Z_{k2^{-n}})+f(Z_{(k+1)2^{-n}})\big)\big[( \Delta_{k,n}Z )^p - 
E[(\Delta_{k,n}Z)^p]\big], \label{N}
\end{equation}
for some $\kappa$ to be discovered. Due to the fact that one cannot separate $X$ from $Y$ inside $Z$ in the definition
of $N_n^{(p)}$, working directly with (\ref{N}) seems to be a  difficult task (see also \cite[ Problem 5.1]{kh-lewis2}). That is why,
following an idea introduced by Khoshnevisan and Lewis \cite{kh-lewis1} in the study of the case $H=1/2$,
we introduce the following collection of stopping times (with
respect to the natural filtration of $Y$), denoted by
\begin{equation}
\mathscr{T}_n=\{T_{k,n}: k\geq 0\}, \quad n\geq 0, \label{T-N}
\end{equation}
which are in turn expressed in terms of the subsequent hitting
times of a dyadic grid cast on the real axis. More precisely, let
$\mathscr{D}_n= \{j2^{-n/2}:\,j\in\Z\}$, $n\geq 0$, be the dyadic
partition (of $\R$) of order $n/2$. For every $n\geq 0$, the
stopping times $T_{k,n}$, appearing in (\ref{T-N}), are given by
the following recursive definition: $T_{0,n}= 0$, and
\[
T_{k,n}= \inf\big\{s>T_{k-1,n}:\quad
Y(s)\in\mathscr{D}_n\setminus\{Y_{T_{k-1,n}}\}\big\},\quad k\geq 1.
\]
  As shown in
\cite{kh-lewis1}, as $n$ tends to
infinity the collection $\{T_{k,n}:\,1\leq k \leq 2^nt\}$ approximates the
common dyadic partition $\{k2^{-n}:\,1\leq k \leq 2^nt\}$ of order $n$ of the time interval $[0,t]$ (see
\cite[Lemma 2.2]{kh-lewis1} for a precise statement).
Based on this fact, one can introduce the counterpart of (\ref{N}) based on $\mathscr{T}_n$, namely,
\begin{equation*}
\tilde{N}_n^{(p)}(t)=2^{-n\tilde{\kappa}}\sum_{k=0}^{\lfloor 2^n t\rfloor-1}\frac12\big(f(Z_{T_{k,n}})+f(Z_{T_{k+1,n}})\big)\big[\big(2^{\frac{nH}{2}}(Z_{T_{k+1,n}}- Z_{T_{k,n}})\big)^p - \mu_p\big], 
\end{equation*}
 with $\mu_p:= E[N^p]$, where $N \sim \mathscr{N}(0,1)$ and for some $\tilde{\kappa}>0$ to be discovered. 
At this stage, it is worthwhile noting that we are dealing with a modified weighted $p$-variation of $Z$. In fact, the collection of stopping times $\{T_{k,n}:\,1\leq k \leq 2^nt\}$  will play an important role in our analysis as we will see in Lemma \ref{lemme-algebric}. 

\subsubsection{Known results about the weighted power variations of fBmBt}
The asymptotic behavior of $\tilde{N}_n^{(p)}(t)$, as $n$ tends to infinity, has been studied in \cite{NP} when $H=1/2$. For $H=1/2$, one can deduce the following  finite dimensional distributions (f.d.d.) convergence in law from \cite[Theorem 1.2]{NP}.
\begin{itemize}
\item[1)] For $f\in C_b^2(\R)$ and for any integer $r \geq 2$,  we have
\end{itemize}
\begin{eqnarray}
&&\bigg( 2^{-\frac{n}{4}}\sum_{k=0}^{\lfloor 2^n t \rfloor -1} \frac12\big(f(Z_{T_{k,n}})+f(Z_{T_{k+1,n}})\big)\big(2^{\frac{n}{4}}(Z_{T_{k+1,n}}-Z_{T_{k,n}})\big)^{2r-1}  \bigg)_{t \geq 0} \notag\\
&&\underset{n\to \infty}{\overset{\rm f.d.d.}{\longrightarrow}} \bigg( \int_0^{Y_t} f(X_s)(\mu_{2r} d^{\circ}X_s + \sqrt{\mu_{4r-2} - \mu_{2r}^2}\,dW_s) \bigg)_{t \geq 0}, \label{NP2}
\end{eqnarray}
with $\mu_n:= E[N^n]$, where $N \sim \mathscr{N}(0,1)$,  for all $t \in \R$, $\int_0^t f(X_s)d^{\circ}X_s$ is the Stratonovich integral of $f(X)$ with respect to $X$ defined as the limit in probability of $2^{-\frac{nH}{2}}W_{n}^{(1)}(f,t)$ as $n \to \infty$, with $W_{n}^{(1)}(f,t)$ defined in (\ref{Wn2}), $W$ is a two-sided Brownian motion independent of $(X,Y)$ and for $u \in \R$, $\int_0^{u} f(X_s)dW_s$ is the Wiener-It\^o integral of $f(X)$ with respect to $W$ defined in (\ref{integrale}).

For $H\neq 1/2$, the second author of this paper has proved in \cite{RZ4} the following result with  $f\in C_b^\infty(\R)$ ($f$  is infinitely differentiable with bounded derivatives of all orders), 
\begin{itemize}
\item[2)] For $\frac16 < H < \frac12$ and for any integer $r\geq 2$, we have
\end{itemize}
\begin{eqnarray}
&&\bigg( 2^{-\frac{n}{4}}\sum_{k=0}^{\lfloor 2^n t \rfloor -1} \frac12\big(f(Z_{T_{k,n}})+f(Z_{T_{k+1,n}})\big)\big(2^{\frac{nH}{2}}(Z_{T_{k+1,n}}-Z_{T_{k,n}})\big)^{2r-1}  \bigg)_{t \geq 0} \notag\\
&&\underset{n\to \infty}{\overset{\rm f.d.d.}{\longrightarrow}} \bigg( \beta_{2r-1}\int_0^{Y_t} f(X_s)dW_s \bigg)_{t\geq 0}, \label{fdd 1}
\end{eqnarray}
where  for $u \in \R$, $\int_0^{u} f(X_s)dW_s$ is the Wiener-It\^o integral of $f(X)$ with respect to $W$ defined in (\ref{integrale}) and $\beta_{2r-1} = \sigma_r$, where $\sigma_r$ is defined in Theorem \ref{intro-main}.
\begin{itemize}
\item[3)]  Fix a time $t \geq 0$, for $H>\frac12$ and for any integer $r \geq 1$, we have
\end{itemize}
\begin{eqnarray}
2^{-\frac{nH}{2}}\sum_{k=0}^{\lfloor 2^n t \rfloor -1} \frac12\big(f(Z_{T_{k,n}})+f(Z_{T_{k+1,n}})\big)\big(2^{\frac{nH}{2}}(Z_{T_{k+1,n}}-Z_{T_{k,n}})\big)^{2r-1} \underset{n\to \infty}{\overset{L^2}{\longrightarrow}}  \frac{(2r)!}{r!2^r}\int_0^{Y_t}f(X_s)d^{\circ}X_s,\notag\\ \label{L2}
\end{eqnarray}
where for all $t \in \R$, $\int_0^t f(X_s)d^{\circ}X_s$ is defined as in (\ref{NP2}). 

As it has been mentioned in \cite{RZ4}, the limit of the  weighted $(2r-1)$-variation of $Z$  for $H=\frac12$ in (\ref{NP2}) is intermediate between the limit of the  weighted $(2r-1)$-variation of $Z$ for $\frac16<H<\frac12$ in (\ref{fdd 1}) and the limit of the  weighted $(2r-1)$-variation of $Z$ for $H> \frac12$ in (\ref{L2}). 
A natural question is then to discovered what happens for $H\leq 1/6$. The answer is given in Theorem \ref{application theorem}.

\begin{remark}
One  can remark that, thanks to Theorem \ref{application theorem},   (\ref{fdd 1}) holds true for $H\leq 1/6$. 
\end{remark}

\subsubsection{Asymptotic behavior of the trapezoidal  weighted odd-power variations of  the fBmBt for $H<1/2$}

The asymptotic behavior of the trapezoidal weighted odd-power variations of  the fBmBt for $H<1/2$ is given in Theorem \ref{application theorem}.  Inspired by \cite{kh-lewis1}, the proof of Theorem \ref{application theorem}, given below,   will be done in several steps.

\subsubsection*{Step 1: A key lemma}
For each
integer $n\geq 1$, $k\in\Z$ and real number $t\geq 0$, let $U_{j,n}(t)$ (resp.
$D_{j,n}(t)$) denote the number of \textit{upcrossings} (resp.
\textit{downcrossings}) of the interval
$[j2^{-n/2},(j+1)2^{-n/2}]$ within the first $\lfloor 2^n
t\rfloor$ steps of the random walk $\{Y_{T_{k,n}}\}_{k\geq 0}$, that is,
\begin{eqnarray}
U_{j,n}(t)=\sharp\big\{k=0,\ldots,\lfloor 2^nt\rfloor -1 :&&
\notag
\\ Y_{T_{k,n}}\!\!\!\!&=&\!\!\!\!j2^{-n/2}\mbox{ and }Y_{T_{k+1,n}}=(j+1)2^{-n/2}
\big\}; \notag\\
D_{j,n}(t)=\sharp\big\{k=0,\ldots,\lfloor 2^nt\rfloor -1:&&
\notag
\\ Y_{T_{k,n}}\!\!\!\!&=&\!\!\!\!(j+1)2^{-n/2}\mbox{ and }Y_{T_{k+1,n}}=j2^{-n/2}
\big\}.\notag
\end{eqnarray}
The following lemma taken from \cite[Lemma 2.4]{kh-lewis1} is going to be the key when
studying the asymptotic behavior of the weighted power variation $V_n^{(2r-1)}(f,t)$ of order $r\geq 1$, defined as:
\begin{equation}
V_n^{(2r-1)}(f,t)=\sum_{k=0}^{\lfloor 2^n t \rfloor -1} \frac12\big(f(Z_{T_{k,n}})+f(Z_{T_{k+1,n}})\big)\big[\big(2^{\frac{nH}{2}}(Z_{T_{k+1,n}}-Z_{T_{k,n}})\big)^{2r-1} \big],\quad t\geq 0. \label{Vn}
\end{equation}
Its main feature is to separate $X$ from $Y$, thus providing a representation of
$V_n^{(2r-1)}(f,t)$ which is amenable to analysis.

\begin{lemme}\label{lemme-algebric}
Fix $f:\R\rightarrow \R$, $t\geq 0$ and $r\in \N^{*}$.
Then
\begin{eqnarray}
&&V_n^{(2r-1)}(f,t) = \sum_{j\in\Z}
\frac12\left(
f(X_{j2^{-\frac{n}{2}}}) + f(X_{(j+1)2^{-\frac{n}{2}}})\right) \big[\big(2^{\frac{nH}{2}}(X_{(j+1)2^{-\frac{n}2}}-
X_{j2^{-\frac{n}{2}}})\big)^{2r-1}  \big] \notag\\
&& \hspace{6cm}\times\big(U_{j,n}(t)-D_{j,n}(t)\big).
\end{eqnarray}
\end{lemme}

\subsubsection*{Step 2: Transforming the weighted power variations of odd order}

By \cite[Lemma 2.5]{kh-lewis1}, one has
\[
U_{j,n}(t) - D_{j,n}(t)= \left\{
\begin{array}{lcl}
 {\mathbf 1}_{\{0\leq j< j^*(n,t)\}} & &\mbox{if $j^*(n,t) > 0$}\\
 0 & &\mbox{if $j^{*}(n,t) = 0$}\\
  -{\mathbf 1}_{\{j^*(n,t)\leq j<0\}} & &\mbox{if $j^*(n,t)< 0$}
  \end{array}
\right.,
\]
where $j^*(n,t)=2^{n/2}Y_{T_{\lfloor 2^n t\rfloor,n}}$.
As a consequence, $V_n^{(2r-1)}(f,t)$ is equal to
\begin{eqnarray*}
\left\{
\begin{array}{lcl}
\sum_{j=0}^{j^*(n,t)-1}\frac12\big(f(X^+_{j2^{-n/2}}) + f(X^+_{(j+1)2^{-n/2}})\big) \big(X^{n,+}_{j+1}- X^{n,+}_j\big)^{2r-1} & &\mbox{if $j^*(n,t) > 0$}\\
  0 &&\mbox{if $j^{*}(n,t) = 0$}\\
\sum_{j=0}^{|j^*(n,t)|-1} \frac12\big(f(X^-_{j2^{-n/2}}) + f(X^-_{(j+1)2^{-n/2}})\big) \big(X^{n,-}_{j+1}- X^{n,-}_j\big)^{2r-1} & &\mbox{if $j^*(n,t) < 0$}
  \end{array}
\right.,
\end{eqnarray*}
where $X^+_t := X_t$ for $t\geq 0$, $X^-_{-t} :=X_t$ for $t<0$, $X^{n,+}_{t} := 2^{\frac{nH}{2}}X^+_{2^{-\frac{n}{2}}t}$ for $ t \geq 0$ and $X^{n,-}_{-t} := 2^{\frac{nH}{2}}X^-_{2^{-\frac{n}{2}}(-t)}$ for $t < 0$.

Let us now
introduce the following sequence of processes $W_{\pm,n}^{(2r-1)}$:
\begin{equation}
W_{\pm,n}^{(2r-1)}(f,t)=
 \sum_{j=0}^{\lfloor 2^{n/2}t\rfloor -1} \frac12\big(f(X^\pm_{j2^{-\frac{n}{2}}}) + f(X^\pm_{(j+1)2^{-\frac{n}{2}}})\big) (X^{n,\pm}_{j+1}- X^{n,\pm}_{j})^{2r-1},\quad t \geq 0 \label{Wn}
 \end{equation}
 \begin{equation}
 W_{n}^{(2r-1)}(f,t):=\left \{ \begin{array}{lc}
                      W_{+,n}^{(2r-1)}(f,t) &\text{if $t \geq 0$}\label{Wn2}\\
                      W_{-,n}^{(2r-1)}(f,-t) &\text{if $t < 0$}
                      \end{array}
                      \right. .
\end{equation}
We then have, 
\begin{eqnarray}
V_n^{(2r-1)}(f,t)
= W_{n}^{(2r-1)}(f,Y_{T_{\lfloor 2^n t\rfloor,n}}).
\label{transforming}
\end{eqnarray}

\subsubsection*{Step 3: A result  concerning the trapezoidal weighted odd-power variations of the fBm}

We have the following  proposition.
 
\begin{prop}\label{proposition principale}  Let  $H < \frac12$. Given an integer $r\geq 1$  then,
for any $f\in C^M(\mathbb{R})$, where $M >2r-2 + \frac 1{2H} $, 
\begin{equation}\label{conv} 
\bigg(2^{-\frac{n}{4}}W_{n}^{(2r-1)}(f, t)\bigg)_{t\in \R} \underset{n\to \infty}{\overset{\rm Law}{\longrightarrow}} \bigg(\sigma_r\int_0^{t}f(X_s)dW_s \bigg)_{t\in \R},
\end{equation}
in $D(\R)$, where $W_{n}^{(2r-1)}(f,t)$ is defined in (\ref{Wn2}), $W$ is a two-sided Brownian motion independent of $(X,Y)$, and  $\int_0^{t} f(X_s)dW_s$ is defined in the following natural way: for $u\in \R$, 
\begin{eqnarray}\label{integrale}
\int_0^uf(X_s)dW_s :=\left\{
\begin{array}{lcl}
\int_0^uf(X^+_s)dW^+_s  &&\mbox{if $ u \geq 0$}\\
 \int_0^{-u}f(X^-_s)dW^-_s  &&\mbox{if $u < 0$}
  \end{array}
\right.,
\end{eqnarray}
where   $W^{+}_t=W_t$ if $t>0$ and
$W^{-}_t=W_{-t}$ if $t<0$,  $X^+$ and $X^-$ are defined in Step 2, and $\int_0^u f(X^{\pm}_s)dW^{\pm}_s$ must be  understood in the Wiener-It\^o sense.
\end{prop}
 
 \noindent
 {\it Proof.}
We define, for all $j,n \in \N$, $\tilde{\beta}^\pm_{j,n}:= \frac12(X^\pm_{j2^{-\frac{n}{2}}}+X^\pm_{(j+1)2^{-\frac{n}{2}}})$.  Let us 
introduce the following sequence of processes:
\begin{equation*}
M_{\pm,n}(f,t)=
 \sum_{j=0}^{\lfloor 2^{n/2}t\rfloor -1} f\big( \tilde{\beta}^\pm_{j,n}\big) (X^{n,\pm}_{j+1}- X^{n,\pm}_{j})^{2r-1},\quad t \geq 0, \label{Mn}
 \end{equation*}
 \begin{equation}
 M_{n}(f,t):=\left \{ \begin{array}{lc}
                      M_{+,n}(f,t) &\text{if $t \geq 0$}\label{Mn2}\\
                      M_{-,n}(f,-t) &\text{if $t < 0$}
                      \end{array}
                      \right. .
\end{equation}
 Then, by the same arguments that have been used in the proof of Proposition \ref{Prop}, we have
 \[
 2^{-\frac{n}{4}}M_{n}(f,\cdot) - 2^{-\frac{n}{4}}W_{n}^{(2r-1)}(f,\cdot)\underset{n \to + \infty}{\longrightarrow} 0,
 \]
in probability in $D(\R)$. So, in order to prove (\ref{conv}) it is enough to prove the following result
\begin{equation}\label{conv-bis} 
\bigg(2^{-\frac{n}{4}}M_{n}(f,t)\bigg) _{t\in \R}\underset{n\to \infty}{\overset{\rm Law}{\longrightarrow}} \bigg(\sigma_r\int_0^{t}f(X_s)dW_s \bigg)_{t \in \R}, 
\end{equation}
in $ D(\R)$.
 The proof of (\ref{conv-bis}) will be done in two steps, first we prove the  convergence in law of the finite dimensional distributions and later  we prove  tightness.
 
 \medskip
\noindent
{\bf  1. Convergence in law of the finite dimensional distributions.}
 Our purpose is to prove that
 \[
\bigg(2^{-\frac{n}{4}}M_{n}(f,t)\bigg)_{ t \in \R} \underset{n\to \infty}{\overset{f.d.d.}{\longrightarrow}} \bigg(\sigma_r\int_0^{t}f(X_s)dW_s \bigg)_{ t \in \R},
 \]
 which is equivalent, by (\ref{Mn2}), to prove that 
 \begin{equation}  \label{conv'} 
\bigg(2^{-\frac{n}{4}}M_{\pm,n}(f,t)\bigg)_{ t \geq 0} \underset{n\to \infty}{\overset{f.d.d.}{\longrightarrow}} \bigg(\sigma_r\int_0^{t}f(X^\pm_s)dW^\pm_s \bigg)_{ t \geq 0}.
\end{equation}
The proof of (\ref{conv'}) uses arguments similar to those  employed in   part {\bf (A)}  of the proof of Theorem \ref{intro-main}, the main ingredient being the  small blocks/big blocks approach. Fix $m\leq n$  and for each $j \geq 0$ we denote by $k:= k(j) =
\sup\{i\geq 0: i2^{-m/2}  \leq j2^{-n/2} \}$.
Define    
\bea
\widetilde{M}^\pm _{n,m} (f,t) =  \sum_{j=0}^{\lfloor 2^{n/2} t\rfloor -1 }  f(\tilde{\beta}^\pm_{k(j),m}) \big( X^{n,\pm}_{j+1}- X^{n,\pm}_{j}\big)^{2r-1}. \notag
\eea
 It is known that (see (3.5) in \cite{RZ1} and  part {\bf (a)} in the proof of Proposition 5.1 in \cite{RZ4}) 
 \begin{equation*} 
\bigg( 2^{-n/4}\widetilde{M}^\pm _{n,m} (f,t)\bigg)_{ t \geq 0} \underset{n\to \infty}{\overset{f.d.d.}{\longrightarrow}} \bigg( L^\pm_{ m}(t) \bigg)_{ t \geq 0},
\end{equation*} 
 where
 \begin{eqnarray*}
 L^\pm_{ m}(t)&:=& \sigma_r  \sum_{k=0}^{\lfloor 2^{m/2} t\rfloor -1} f(\tilde{\beta}^\pm_{k,m})\big(  W^\pm_{(k+1)2^{-m/2} }- W^\pm_{k2^{-m/2}}\big) \\
&& + \sigma_r  f(\tilde{\beta}^\pm_{\lfloor 2^{m/2} t \rfloor ,m})
 \big (W^\pm_{t}- W^\pm_{(\lfloor 2^{m/2} t \rfloor ) 2^{-{m/2}} }\big),
\end{eqnarray*}
 with $W^{+}_t=W_t$ if $t>0$ and
$W^{-}_t=W_{-t}$ if $t<0$, where $W$ is a two-sided Brownian motion independent of $(X,Y)$. From the theory of stochastic calculus for semimartingales, we deduce that $L^\pm_{ m}(t) \overset{L^2}{\longrightarrow}\sigma_r \int_0^{t}f(X^\pm_s)dW^\pm_s $ as $m \to \infty$. Then, it is enough to prove that, for all $t\geq 0$,
\begin{align*}
\lim_{m\rightarrow\infty}\limsup_{n\rightarrow\infty}\| 2^{-\frac{n}{4}}M_{+,n}(f,t)-2^{-n/4}\widetilde{M}^+ _{n,m} (f,t) \|_{L^{2}(\Omega)}
  &=0,\\
  \lim_{m\rightarrow\infty}\limsup_{n\rightarrow\infty}\| 2^{-\frac{n}{4}}M_{-,n}(f,t)-2^{-n/4}\widetilde{M}^- _{n,m} (f,t) \|_{L^{2}(\Omega)}
  &=0.
\end{align*}
The proof of the last claim is similar to the proof of (\ref{lim:PhiPhitilde}) and is left to the reader.

 \medskip
 \noindent
{\bf 2. Proof of  Tightness.}
We claim that the distribution of the sequence  $\big(2^{-\frac{n}{4}}M_{n}(f,\cdot)\big)_{n \in \N}$ is tight in $D(\R)$. 
To prove this claim we  will show   that for any $T >0$ and for every $-T< s\leq t < T$, and $p>2$, there exists a constant $C>0$, such that
\bea \label{claim'}
E\left[ |2^{-\frac{n}{4}}M_{n}(f,t)-2^{-\frac{n}{4}}M_{n}(f,s)|^{p}\right]\leq C \left(\frac{\floor{2^n t}-\floor{2^n s}}{2^n}\right)^{\frac{p}{2}}
+ C\left(\frac{\floor{2^n t}-\floor{2^n s}}{2^n}\right)^{pH}.
\eea
To do so, we  distinguish three cases, according to the sign of $s, t \in \R$:

\medskip
\noindent
{\bf i)}. Suppose that $0 \leq s \leq t$. In this case we can write
\begin{eqnarray*}
 && E\left[ |2^{-\frac{n}{4}}M_{n}(f,t)-2^{-\frac{n}{4}}M_{n}(f,s)|^{p}\right]\\
&=& E\left[ |2^{-\frac{n}{4}}M_{+,n}(f,t)-2^{-\frac{n}{4}}M_{+,n}(f,s)|^{p}\right]\\
&\leq & C \left(\frac{\floor{2^n t}-\floor{2^n s}}{2^n}\right)^{\frac{p}{2}}
+ C\left(\frac{\floor{2^n t}-\floor{2^n s}}{2^n}\right)^{pH},
\end{eqnarray*}
where the proof of the last inequality is the same  as the proof of (\ref{claim}).

\medskip
\noindent
{\bf ii)}.
 Suppose $s \leq t \leq 0$. Then, we have
\begin{eqnarray*}
&&E\left[ |2^{-\frac{n}{4}}M_{n}(f,t)-2^{-\frac{n}{4}}M_{n}(f,s)|^{p}\right]\\
&=&E \left[ |2^{-\frac{n}{4}}M_{-,n}(f,-t)-2^{-\frac{n}{4}}M_{-,n}(f,-s)|^{p}\right]\\
&\leq & C \left(\frac{\floor{2^n (-s)} - \floor{2^n (-t)}}{2^n}\right)^{\frac{p}{2}}
+ C\left(\frac{\floor{2^n (-s)} - \floor{2^n (-t)}}{2^n}\right)^{pH}\\
&=& C \left(\frac{\floor{2^n t}-\floor{2^n s}}{2^n}\right)^{\frac{p}{2}}
+ C\left(\frac{\floor{2^n t}-\floor{2^n s}}{2^n}\right)^{pH},
\end{eqnarray*}
where the proof of the second  inequality is the same as the  proof of (\ref{claim}) and we get the last equality  since for any $x <0$, $\lfloor -x \rfloor = - \lfloor x \rfloor -1$.

\medskip
\noindent
{\bf iii)}.  Suppose $ s < 0 < t$. Then, we can write
\begin{eqnarray*}
&& E\left[ |2^{-\frac{n}{4}}M_{n}(f,t)-2^{-\frac{n}{4}}M_{n}(f,s)|^{p}\right]\leq C\big( E\left[ |2^{-\frac{n}{4}}M_{n}(f,t)-2^{-\frac{n}{4}}M_{n}(f,0)|^{p}\right] \\
&&+ E\left[ |2^{-\frac{n}{4}}M_{n}(f,s)-2^{-\frac{n}{4}}M_{n}(f,0)|^{p}\right]\big)\\
&=&C\big( E\left[ |2^{-\frac{n}{4}}M_{+,n}(f,t)-2^{-\frac{n}{4}}M_{+,n}(f,0)|^{p}\right]\\
&& +E\left[ |2^{-\frac{n}{4}}M_{-,n}(f,-s)-2^{-\frac{n}{4}}M_{-,n}(f,0)|^{p}\right]\big)\\
&\leq & C \left(\frac{\floor{2^n t}}{2^n}\right)^{\frac{p}{2}}
+ C\left(\frac{\floor{2^n (-s)}}{2^n}\right)^{pH} \\
&\leq& C \left(\frac{\floor{2^n t}+ \floor{2^n (-s)} +1}{2^n}\right)^{\frac{p}{2}}
+ C\left(\frac{\floor{2^n (-s)} + \floor{2^n t}+1}{2^n}\right)^{pH}\\
&=&C \left(\frac{\floor{2^n t}-\floor{2^n s}}{2^n}\right)^{\frac{p}{2}}
+ C\left(\frac{\floor{2^n t}-\floor{2^n s}}{2^n}\right)^{pH},
\end{eqnarray*}
where we have the  third inequality by  i) and ii).

Finally, we have proved (\ref{claim'}) which proves the tightness of $\big(2^{-\frac{n}{4}}M_{n}(f,\cdot)\big)_{n \in \N}$ in $D(\R)$.

\subsubsection*{Step 4: Convergence in law of $Y_{T_{\lfloor 2^n \cdot \rfloor,n}}$}
As it has been mentioned in \cite{kh-lewis1},   $\{2^{n/2}Y_{T_{k,n}} : k \geq 0\}$ is a simple and symmetric random walk on  $\Z$. Observe that for all $t\geq 0$, $Y_{T_{\lfloor 2^n t\rfloor,n}} = 2^{-n/2}\times 2^{n/2}Y_{T_{\lfloor 2^n t\rfloor,n}} = 2^{-n/2} \sum_{l=0}^{\lfloor 2^n t\rfloor -1} 2^{n/2}(Y_{T_{l+1,n}} - Y_{T_{l,n}})$, where $\big(2^{n/2}(Y_{T_{l+1,n}} - Y_{T_{l,n}})\big)_{l \in \N}$ are independent and identically distributed random variables  following the Rademacher distribution. By Donsker theorem, we get that

\begin{equation}
\big(Y_{T_{\lfloor 2^n t \rfloor,n}}\big)_{t \geq 0} \underset{n \to  \infty}{\overset{law}{\longrightarrow}} (Y_t)_{t \geq 0} \text{\: in \:} D([0, +\infty)). \label{conv-law}
\end{equation}

\subsubsection*{Step 5: Last step in the proof of Theorem \ref{application theorem}}
 Thanks to Proposition \ref{proposition principale}, to (\ref{conv-law}), and to the independence of $X$, $W$ and $Y$ , we have
 \begin{equation}\label{last formula}
 \big(2^{-\frac{n}{4}}W_{n}^{(2r-1)}(f,\cdot), Y_{T_{\lfloor 2^n \cdot \rfloor,n}}\big) \underset{n \to +\infty}{\overset{law}{\longrightarrow}} \big(\sigma_r\int_0^{\cdot}f(X_s)dW_s  , Y) \text{\: in \:} D(\R)\times D([0, +\infty)). 
 \end{equation}
 Let us define $(B_t)_{t \in \R}$ as follows $B_t:=\sigma_r\int_0^{t}f(X_s)dW_s$. 
Since $(x,y) \in D(\R)\times D([0, +\infty)) \mapsto x\circ y \in D([0, +\infty))$ is measurable (see M16 at page 249 in \cite{Billingsley} for a proof of this result) and since $B \circ Y$ is continuous, then, by (\ref{last formula}) and Theorem 2.7  in \cite{Billingsley}, it follows that
\[
2^{-\frac{n}{4}}W_{n}^{(2r-1)}(f,Y_{T_{\lfloor 2^n \cdot \rfloor,n}}) \underset{n \to +\infty}{\overset{law}{\longrightarrow}} B \circ Y = \sigma_r\int_0^{Y(\cdot)}f(X_s)dW_s, \text{\; in \:} D([0, +\infty)).
\]
The proof of Theorem \ref{application theorem} follows from (\ref{transforming}) and the last convergence in law.

\bigskip
\noindent
{\bf Acknowledgment}.  The first author  was supported by the NSF grant DMS 1512891. The first drafted version of this paper has been done when the second author was a member of the Research training group 2131, working in the Technical University of Dortmund. He is  thankful to the financial support of the DFG (German Science Foundations) Research Training Group 2131. He is also thankful to the financial support of FRIAS-USIAS program.

\end{document}